%% file: IntervalsAndOuterMeasureOnR.tex
\documentclass[a4paper,portrait]{article}

\usepackage[a4paper, top=0.5in, bottom=0.8in, hmargin=1.2in, footnotesep=0.55in]{geometry}
\usepackage{amsmath}
\usepackage{amssymb}
\usepackage{tikz}
\usepackage{enumitem}
\usepackage[title, toc, titletoc]{appendix}

\usepackage{mathdots}

\usepackage{hyperref}
\hypersetup{
	hypertexnames=false, 
	colorlinks=true,
	allcolors=black,
    urlcolor=blue!50!black,
	bookmarksnumbered=true,
	bookmarksopen=true,
	pdfauthor=Ross Ure Anderson
}

\usetikzlibrary{intersections}
\usetikzlibrary{calc}
\usetikzlibrary{positioning}
\usetikzlibrary{arrows.meta}

\newcommand{\tf}{\raisebox{0.2ex}{$\therefore $}}
\newcommand{\imp}{\ensuremath{\Rightarrow}}
\newcommand{\ra}{\ensuremath{\rightarrow}}
\newcommand{\rn}{\ensuremath{\mathbb{R}}}
\newcommand{\ern}{\ensuremath{\overline{\mathbb{R}}}}
\newcommand{\nn}{\ensuremath{\mathbb{N}}}
\newcommand{\ci}{\ensuremath{\mathcal{I}}}
\newcommand{\cb}{\ensuremath{\mathcal{B}}}
\newcommand{\cc}{\ensuremath{\mathcal{C}}}
\newcommand{\ca}{\ensuremath{\mathcal{A}}}
\newcommand{\cl}{\ensuremath{\mathcal{L}}}
\newcommand{\sg}{\ensuremath{\sigma}}
\newcommand{\ep}{\ensuremath{\epsilon}}
\newcommand{\sa}{$\sigma$-algebra}
\newcommand{\sm}{\setminus}
\newcommand{\es}{\O}
\newcommand{\leftcurly}[1]{\left\{ \rule{0ex}{#1} \right.}
\newcommand{\rightcurly}{\left. \rule{0ex}{3ex} \right\}}
\newcommand{\leftvert}[1]{\left| \rule{0ex}{#1} \right.}
\newcommand{\rightvert}[1]{\left. \rule{0ex}{#1} \right|}

\newenvironment{proof}{{\noindent \large \bfseries Proof} \par \vspace{1.3ex}}{\par \hfill QED \vspace{2.7ex}}
\newenvironment{proof-sm}{{\noindent \bfseries Proof} \par \vspace{1.3ex}}{\par \hfill QED \vspace{2.7ex}}
\newenvironment{note}{{\noindent \large \bfseries Note} \par \vspace{1.3ex}}{\par}
\newenvironment{notes}{{\noindent \large \bfseries Notes} \par \vspace{0ex}}{\par}
\newenvironment{corollary}{{\noindent \large \bfseries Corollary} \par \vspace{2ex}}{\par}

\newcommand{\terminatingbox}{\noindent \hfill \raisebox{0ex}{\framebox[0.6em]{\rule{0em}{0.6ex}}}}

\title{\vspace{-2ex}Intervals and Outer Measure on \rn}
\author{Ross Ure Anderson \thanks{E-mail: ruanderson100@yahoo.com; web: \url{https://archive.org/details/@ross_ure_anderson}}}
\date{2\textsuperscript{nd} October, 2023}

\begin{document}

\maketitle

\section{Introduction}

This article gives some properties of intervals in \rn\ and discusses some problems involving intervals for which the concept of outer measure on \rn\ provides a more efficient solution than an elementary approach. The outer measure is then defined and some of its main properties in relation to intervals are developed, culminating in the countable additivity of outer measure on the `system of intervals' \ci\ $= \{$ all countable unions of intervals in \rn\ $\}$. This demonstrates early on how the outer measure on \rn\ is naturally countably additive on a quite large class of sets, and motivates the Borel algebra \cb\ as an extension of that class which provides an additional desired property of outer measure, namely closure of its domain under set complementation --- for example as developed in \cite[Chap~2]{axler:measure-theory}. Details are given of how one of the intervals problems solved by the outer measure allows proof prior to the Lebesgue integration theory of the Bounded Convergence, Monotone Convergence, and Dominated Convergence Theorems for Riemann integrals. One application of the latter is the proof of Stirling's Formula in \cite{conrad:stirlings-formula}. Some further details on handling double series are provided than is normally given, based on the textbook \cite{knopp:infinite-series} of Konrad Knopp and the article \cite{mse:series-rearrangements}. The term `countable union' of sets will mean a union of an infinite sequence of sets. \ern\ will denote the extended real number system $\rn\ \cup \{\infty, -\infty\}$.

\section{Properties of Intervals}

\subsection{Definition (Interval Length)}
\label{defn:interval-length}
The length $l(I) \in [0, \infty]$ of an interval $I$ in \rn\ is given by :
$$
l(I) = \left\{
\begin{array}{ll}
b - a & \mbox{if } I = (a, b), [a, b], [a, b), or (a, b] \mbox{ where } a, b \in \rn \mbox{ and } a \leq b, \\
\infty & \mbox{if } I = (a, \infty), [a, \infty), (-\infty, a), \mbox{ or } (-\infty, a] \mbox{ where } a \in \rn, \\
\infty & \mbox{if } I = (-\infty, \infty), \\
0 & \mbox{if } I = \es
\end{array}
\right.
$$
\noindent Thus from the arithmetic of \ern\ :
$$l(I) = \sup I - \inf I, \hspace{2em} \forall \mbox{ non-empty intervals } I \subseteq \rn$$
and $l(I) = \infty$ iff $I$ is unbounded.

\subsection[Interval Theorem 2.2]{Theorem}
\label{thm:propints-finite-disj}
If $I_1, \ldots, I_n$ are disjoint subintervals of interval $I$ then :
$$l(I_1) + \cdots + l(I_n) \leq l(I).$$

\begin{proof}
\noindent The result is clear if $l(I) = \infty$. Otherwise $l(I) < \infty$ and $I$ is a bounded interval. Relabel the $I_1, \ldots, I_n$ so that they are in left to right order. Then whatever types of bounded intervals $I_k$ are (closed, open, half-closed/half-open), if $g \geq 0$ is the total gap length within $I$ between the $I_k$ then~:
$$l(I_1) + \ldots + l(I_n) + g = l(I)$$
and thus the required result follows.
\end{proof}

\subsection[Interval Theorem 2.3]{Theorem}
\label{thm:propints-infinite-disj}
If $(I_n)$ is a sequence of disjoint subintervals of interval $I$ then :
$$\sum_{n=1}^{\infty} l(I_n) \leq l(I).$$

\begin{proof}
\noindent The result is clear if $l(I) = \infty$. Otherwise $\forall \; n$ we have from Theorem \ref{thm:propints-finite-disj} above :
$$l(I_1) + \cdots + l(I_n) \leq l(I)$$
within \rn, hence taking the limit as $n \ra \infty$ the result follows.
\end{proof}

\vspace{2ex}

\subsection[Interval Theorem 2.4]{Theorem}
\label{thm:propints-subtract-1-int}
If $I$ and $J$ are intervals then $I \sm J$ is either a subinterval of $I$ or a disjoint union of two subintervals of $I$. \\

\begin{proof}
\noindent For all the possible interval types and relative positions of the $I$ and $J$ the result is clear.
\end{proof}

\subsection[Interval Theorem 2.5]{Theorem}
\label{thm:propints-subtract-n-ints}
If $I, J_1, \ldots, J_n$ are intervals then $I \sm (J_1 \cup \cdots \cup J_n)$ is a disjoint union of finitely many subintervals of $I$. \\

\begin{proof}
\noindent True for $n = 1$ by Theorem \ref{thm:propints-subtract-1-int}. Suppose true for $n$ and consider the case of $n + 1$. Then $I \sm (J_1 \cup \cdots \cup J_{n+1}) = K \sm (J_{n+1})$, where $K = I \sm (J_1 \cup \cdots \cup J_n)$ is, by the inductive assumption, a disjoint union of finitely many subintervals of $I$, $K = K_1 \cup \cdots \cup K_m$, say. Then by Theorem \ref{thm:propints-subtract-1-int}, $K \sm J_{n+1} = (K_1 \sm J_{n+1}) \cup \cdots \cup (K_m \sm J_{n+1})$ is a disjoint union of sets each one of which is a union of one or two disjoint subintervals of $I$. But since the $K_r \sm J_{n+1}$ are themselves disjoint so $K \sm J_{n+1}$ is then a disjoint union of subintervals of $I$, as required.
\end{proof}

\subsection[Interval Theorem 2.6]{Theorem}
\label{thm:propints-disj-count-union}
A countable union of intervals is a disjoint countable union of intervals whose sum of lengths does not exceed the original sum of interval lengths. A finite union of intervals is a disjoint finite union of intervals whose sum of lengths does not exceed the original sum of interval lengths. \\

\begin{proof}
\noindent Consider the countable case. Let $A = \bigcup_{n=1}^{\infty} I_n$, where $I_n$ are arbitrary intervals in \rn. Then :
$$
A = \bigcup_{n=1}^{\infty} I_n \sm (I_1 \cup \cdots \cup I_{n-1}), \hspace{2em} \mbox{a disjoint union of sets.}
$$
By Theorem \ref{thm:propints-subtract-n-ints}, each individual set ${I_n \sm (I_1 \cup \cdots \cup I_{n-1})}$ is a disjoint union of finitely many intervals, ${J_1^{(n)} \cup \cdots \cup J_{r_n}^{(n)}}$ say, where $r_n \geq 0$, and where if $J_i^{(n)} \cap J_k^{(m)} \neq \es$ then $n = m$ and $i = k$, thus : \\

\input{diagram-props-intervals}

\vspace{2ex}

and thus taking the diagonal sequence of intervals shown in red, $(K_n)$ say, these are mutually disjoint intervals with $A = \bigcup_{n=1}^{\infty} K_n$, as required. \\

For the sum of interval lengths, this is clear if we have $\sum_{n=1}^{\infty} l(I_n) = \infty$, otherwise $r^{(n)} = {l(J_1^{(n)}) + \cdots + l(J_{r_n}^{(n)})} \leq l(I_n) < \infty\ \forall\ n$, and thus the series $\sum_{n=1}^{\infty} r^{(n)}$ of total row lengths converges in \rn\ to a number $\leq \sum_{n=1}^{\infty} l(I_n)$, so by \cite[Theorem~3]{mse:series-rearrangements} (and adding empty sets onto the end of the rows), the diagonal series of lengths $\sum_{n=1}^{\infty} l(K_n)$ converges to this number also. \\

For the finite case, just replace the $\infty$ in the above argument with an $N \in \nn$ to obtain $A$ equal to a finite disjoint union of intervals, and the sum of interval lengths is clearly not increased.
\end{proof}

\subpdfbookmark{Corollary}{Corollary}
\begin{corollary}
\noindent The system of intervals $\ci = \{ $ all \emph{disjoint} countable unions of intervals in \rn\ $ \}$. 
\end{corollary}

\terminatingbox

\vspace{7ex}

In Theorem \ref{thm:propints-finite-disj} above we considered a finite disjoint union of subintervals ${I_1 \cup \cdots \cup I_n \subseteq I}$. In the next two theorems we consider the corresponding theorems for the cases ${I_1 \cup \cdots \cup I_n = I}$ and ${I_1 \cup \cdots \cup I_n \supseteq I}$ (superset), though disjointness is not required in the latter case.

\subsection[Interval Theorem 2.7]{Theorem}
\label{thm:propints-finite-disj-equal}
If $I_1, \ldots, I_n$ are disjoint subintervals of interval $I$ with ${I = I_1 \cup \cdots \cup I_n}$, then $l(I_1) + \cdots + l(I_n) = l(I)$. \\

\begin{proof}
\noindent Case $l(I) = \infty$, ie. $I$ unbounded. Then at least one $I_k$ is unbounded, so that $l(I_k) = \infty$, \tf\ result follows. Consider the case $l(I) < \infty$, so $I$ and each $I_k$ is bounded. Then $I$ is of form $[a, b], (a, b), [a, b), \mbox{or} (a, b], (a, b \in \rn$), and each subinterval $I_k \subseteq I$ is one of these interval types also. Relabel the $I_k$'s so they are in L to R order. Then in order to avoid a gap the LH of $I_1$ must align with LH of $I$ and be of the same end point type (ie. open or closed). And $\forall \; k$ likewise there can be no gap between RH of $I_k$ and LH of $I_{k+1}$, and the RH of $I_n$ must align with the RH of $I$ and be of the same end point type. Therefore the result follows.
\end{proof}

\newpage
\subsection[Interval Theorem 2.8]{Theorem}
\label{thm:propints-finite-disj-superset}
Suppose interval ${I \subseteq I_1 \cup \cdots \cup I_n}$, a finite (not necessarily disjoint) union of intervals. Then ${l(I_1) + \cdots + l(I_n) \geq l(I)}$. \\

\begin{proof}
\noindent If one of the $I_k$ is unbounded then $l(I_k) = \infty$ and the result is clear, so assume each $I_k$ (and hence $I$) is bounded. \\

Let interval $I'_k = I \cap I_k$ so $I = I'_1 \cup \cdots \cup I'_n$ and therefore :
\begin{eqnarray*}
I & = & I'_1 \cup (I'_2 \sm I'_1) \cup \cdots \cup (I'_n \sm (I'_1 \cup \cdots \cup I'_{n-1}) ) \\
& = & \mbox{union of disjoint sets.}
\end{eqnarray*}

By Theorem \ref{thm:propints-subtract-n-ints}, $\forall \; k, I'_k \sm (I'_1 \cup \cdots \cup I'_{k-1})$ is a disjoint union of finitely many subintervals of $I'_k$, $J_1^{(k)} \cup \cdots \cup J_{r_k}^{(k)}$ say, and by Theorem \ref{thm:propints-finite-disj} we have :
\begin{equation}
\label{eq:propints-superset-1}
l(J_1^{(k)}) + \cdots + l(J_{r_k}^{(k)}) \leq l(I'_k), \hspace{2em} \forall \; k
\end{equation}

Now 
\begin{equation}
\label{eq:propints-superset-2}
I = \bigcup_{k=1}^n J_1^{(k)} \cup \cdots \cup J_{r_k}^{(k)}
\end{equation}
where $J_i^{(k)} \cap J_j^{(l)} \neq \es \imp k = l \mbox{ and } i = j$, so that these sets are mutually disjoint intervals, and so by Theorem \ref{thm:propints-finite-disj-equal} applied to (\ref{eq:propints-superset-2}) we have :
\begin{eqnarray*}
l(I) & = & \sum_{k=1}^n l(J_1^{(k)}) + \cdots + l(J_{r_k}^{k}) \\
& \leq & \sum_{k=1}^n l(I'_k) \hspace{2em} \mbox{from (\ref{eq:propints-superset-1})} \\
& \leq & \sum_{k=1}^n l(I_k) \hspace{2em} \mbox{as } I'_k \subseteq I_k
\end{eqnarray*}
\end{proof}

\begin{notes}
\begin{enumerate}[label=(\arabic*), ref=\arabic*]
\item Theorem \ref{thm:propints-finite-disj-equal} implies finite additivity of the length function on intervals, when applied to intervals equal to disjoint unions of intervals. Theorem \ref{thm:propints-finite-disj-superset} implies finite subadditivity of the length function on intervals, when applied to intervals equal to arbitrary unions of intervals.
\item We find when we develop the properties of the outer measure below that Theorems \ref{thm:propints-finite-disj-equal} and \ref{thm:propints-finite-disj-superset} still apply for $n = \infty$, ie we have (respectively) countable additivity and countable subadditivity of the interval length function. The former follows from the equality of the outer measure to the length function on the intervals (Theorems \ref{thm:om-len-bnd-interval} and \ref{thm:om-len-unbnd-interval}), together with the countable additivity of the outer measure on all intervals (Theorem \ref{thm:om-count-add-all-int}). The latter follows from the equality of the outer measure to the length function on the intervals, together with the countable subadditivity of the outer measure on $P(\rn)$ (Theorem~\ref{thm:om-count-subadd}).
\label{note:propints-note2}
\item Axler~\cite[p20,~\S2.14]{axler:measure-theory} provides an alternative method of proof of Theorem \ref{thm:propints-finite-disj-superset} using induction, in the special case where $I$ is a closed bounded interval and the $I_k$ are arbitrary open intervals.
\end{enumerate}
\end{notes}

\vspace{4ex}

\newpage
\subsection{The Need For Outer Measure}

\label{sec:need-for-outer-measure}

The outer measure provides efficient solutions to the following problems involving intervals, saving effort compared with elementary methods. \\

\noindent \textbf{Example 1} \\

\noindent In Note (\ref{note:propints-note2}) of \S\ref{thm:propints-finite-disj-superset} above we mentioned how outer measure could allow the extension to the infinite case for Theorem \ref{thm:propints-finite-disj-equal} : \\

If $(I_n)$ are disjoint subintervals of $I$ with $I = \bigcup_{n=1}^{\infty} I_n$ then $\sum_{n=1}^{\infty} l(I_n) = l(I)$. \\

It is instructive to consider how we might attempt to prove this without the concept of outer measure. In \cite[Problem~1]{mo:int-problems-qu} a proof is given for the special case where the $I_n$ can be rearranged into L to R order, or into R to L order. In \cite{mo:int-problems-ans1} (a one-dimensional version of \cite[Proposition~2.16,~p48]{weaver:measure-theory}) and \cite{mo:int-problems-ans2} proofs are given for the general case, for $I$ bounded and $I$ unbounded respectively. \\

However using outer measure an immediate proof is obtained by applying Theorems~\ref{thm:om-len-bnd-interval},~\ref{thm:om-len-unbnd-interval}, and \ref{thm:om-count-add-all-int} : simply write $\sum_{n=1}^{\infty} l(I_n) = \sum_{n=1}^{\infty} |I_n| = |\bigcup_{n=1}^{\infty} I_n| = |I| = l(I)$.

\vspace{4ex}

\noindent \textbf{Example 2} \\

\noindent In answer \cite{mo:obvious-hard-to-prove-ans} to question \cite{mo:obvious-hard-to-prove-qu} a similar problem to Example 1 is given --- where intervals $I_n \subseteq [0, 1]\ \forall\ n$ ($I_n$ not necessarily disjoint) and the problem is to show $\sum_{n=1}^{\infty} l(I_n) < 1 \imp \bigcup_{n=1}^{\infty} I_n \subsetneq [0, 1]$. Without the concept of outer measure this is not a trivial statement but it is sometimes mistaken as being `obvious'. An elementary solution which does not use measure theory is provided by the proof in \cite{mo:int-problems-ans1}, which shows that $\sum_{n=1}^{\infty} l(I_n) \geq l(I)$ whenever $I \subseteq \bigcup_{n=1}^{\infty} I_n$ ($I_n$ not necessarily disjoint). \\

With outer measure however we obtain an immediate proof using countable subadditivity (Theorem~\ref{thm:om-count-subadd}) together with extension of the function $l$ (Theorems~\ref{thm:om-len-bnd-interval} and \ref{thm:om-len-unbnd-interval}) : $|\bigcup_{n=1}^{\infty} I_n| \leq \sum_{n=1}^{\infty} |I_n| = \sum_{n=1}^{\infty} l(I_n) < 1 = |[0, 1]|$, so certainly $\bigcup_{n=1}^{\infty} I_n \subsetneq [0, 1]$. (Note: we are thinking of the elementary solution as comprising the text of proof \cite{mo:int-problems-ans1} rather than just using the result of that proof, in comparing it with the outer measure solution where we are assuming the outer measure properties).

\vspace{4ex}

\noindent \textbf{Example 3} \\

\noindent In the proof of the Bounded Convergence Theorem for Riemann Integrals in Bartle \cite[\S22.14,~p288]{bartle:ERA} the following problem arises for which the outer measure provides an efficient solution : prove that if $\delta > 0$ and if $(E_n)$ is a sequence of sets in $[0, 1]$ with each $E_n$ containing a finite number of non-overlapping closed intervals with a total length $\geq \delta$, then there exists a point belonging to infinitely many of the $E_n$. The property of outer measure used to prove this is `measure of a decreasing intersection', as described in \cite[\S2.60,~p44]{axler:measure-theory}. This property applies to outer measure on \rn\ once it is established that outer measure satisfies all the conditions to be a `measure' (\cite[\S2.54,~p41]{axler:measure-theory} and \cite[\S2.23,~p26]{axler:measure-theory}), on a collection of sets in \rn\ containing the `system of intervals' \ci\ (the Borel sets \cb) --- this is carried out in \cite[\S2D]{axler:measure-theory}. This solution is described in \cite[Problem~2]{mo:int-problems-qu}, \cite{mse:countable-additivity}, and Appendix~\ref{sec:convergence-theorems}.

\vspace{4ex}

\noindent \textbf{Example 4} \\

\noindent Show that the `system of intervals' \ci\ is \emph{not} closed under set complementation. The outer measure on \rn\ provides proof of this using the Cantor Set as a counterexample --- see Appendix~\ref{sec:sys-intervals-not-closed}.

\vspace{6ex}

\begin{note}
\noindent In the Examples 1~\&~2 above, all the outer measure properties needed are derived in this article, namely Theorems~\ref{thm:om-count-subadd},~\ref{thm:om-len-bnd-interval},~\ref{thm:om-len-unbnd-interval},~and~\ref{thm:om-count-add-all-int}. \\

For Example~3, the property of outer measure required is `measure of a decreasing intersection' \cite[\S2.60,~p44]{axler:measure-theory} which requires the full development of the outer measure as a `measure' (\cite[\S2.54,~p41]{axler:measure-theory} and \cite[\S2.23,~p26]{axler:measure-theory}) on a suitable class of sets in \rn\ (the Borel sets \cb). The proof of this property requires countable additivity of the outer measure and closure of its domain under set complementation. Outer measure satisfies countable additivity on the set \ci\ by Theorem~\ref{thm:om-count-add-sys-int}, but \ci\ is not closed under set complementation (see Appendix~\ref{sec:sys-intervals-not-closed}) and thus is not a suitable domain for outer measure, but extending \ci\ to \cb\ achieves this closure, whilst maintaining the countable additivity \cite[\S2D]{axler:measure-theory}. Although the outer measure must be restricted to \cb\ to satisfy all the necessary properties of a `measure', note it is well defined on any subset of \rn\ (definition~\S\ref{defn:outer-measure}). \\

Example~4 also makes use of the full properties of outer measure from Axler \cite{axler:measure-theory}.
\end{note}

\section{Outer Measure}

\subsection[Definition (Outer Measure)]{Definition}
\label{defn:outer-measure}
The outer measure $|A| \in [0, \infty]$ of a set $A \subseteq \rn$ is defined by
$$|A| = \inf \leftcurly{3ex} \sum_{k=1}^{\infty} l(I_k) : \{ I_1, I_2, \ldots \} \mbox{ is a countable open interval cover of A} \rightcurly$$

\terminatingbox

\vspace{2ex}

\begin{note}
\noindent The outer measure function is defined on \emph{all} subsets of \rn, in contrast to the length function $l$ above (Definition \ref{defn:interval-length}) which is only defined on intervals in \rn. In Theorems \ref{thm:om-len-bnd-interval} and \ref{thm:om-len-unbnd-interval} below we show the outer measure function is an extension of the function $l$.
\end{note}

\subsection{Theorem (Set Order Preservation)}
\label{thm:om-order-pres}
If $A \subseteq B \subseteq \rn$ then $|A| \leq |B|$. \\

\begin{proof}
\noindent Any countable open interval cover for $B$ is a countable open interval cover for $A$ and hence the sum of its lengths is $\geq |A|$, which is therefore a lower bound for all such sums. But $|B|$ is the greatest such lower bound.
\end{proof}

\subsection{Theorem (Countable Subadditivity)}
\label{thm:om-count-subadd}
For any sequence $(A_n)$ of sets in $\mathbb{R}$ 
\begin{equation}
\leftvert{2.9ex} \bigcup_{n=1}^{\infty} A_n \rightvert{2.9ex} \leq \sum_{n=1}^{\infty} |A_n|
\label{eq:om-count-subadd}
\end{equation}

\begin{proof}
\vspace{2ex}
\noindent Let $A = \bigcup_{n=1}^{\infty} A_n$. In the case where $|A_n| = \infty$ for some $n$ the RHS equals $\infty$, therefore the result is concluded. Assume now that we have every $|A_n| < \infty$. If we have a countable open interval (COI) cover $(I_k^{(n)})$ for $A_n$ ($\forall n$) then : \\

\input{diagram-countable-subadd-1}
{ \vspace{-6.7ex} \begin{equation} \label{eq:om-count-subadd-COIs} \end{equation} } 

\vspace{2ex}

Let $(J_n)$ be the diagonal sequence of open intervals shown in red. Then $(J_n)$ is a COI cover for $A$, and thus by Definition \ref{defn:outer-measure} of outer measure : 
\begin{equation}
|A| \leq \sum_{n=1}^{\infty} l(J_n)
\label{eq:om-count-subadd-diag-sum}
\end{equation}

Intuitively we expect the term on the right of (\ref{eq:om-count-subadd-diag-sum}) to equal the sum of row sums of the interval lengths in the above COI covers (\ref{eq:om-count-subadd-COIs}), and if we choose these COI covers so that for each $n$ the sum of interval lengths in the COI cover for $A_n$ becomes closer and closer to $|A_n|$, the RHS of (\ref{eq:om-count-subadd-diag-sum}) approaches closer and closer to $\sum_{n=1}^{\infty} |A_n|$ from above, thus enforcing the inequality (\ref{eq:om-count-subadd}). \\

To formalize this take an $\epsilon > 0$. Then $\forall \; n$ we can find a COI cover $(I_k^{(n)})$ of $A_n$ such that :
\begin{equation}
|A_n| \leq \sum_{k=1}^{\infty} l(I_k^{(n)}) < |A_n| + \epsilon / 2^n < \infty
\label{eq:om-count-subadd-epsilon-COI}
\end{equation}

Then consider the double series formed by the interval lengths $(l(I_k^{(n)}))$ ($n$ = row, $k$ = column) : \\

\input{diagram-countable-subadd-2}

\vspace{1ex}

where :
\begin{equation}
\mbox{sum of $n$\textsuperscript{th} row } r^{(n)} \in [\;|A_n|, |A_n| + \epsilon / 2^n \;)
\label{eq:om-count-subadd-row-n}
\end{equation}

\vspace{3ex}

From \cite[Theorem~3]{mse:series-rearrangements}, $\sum_{n=1}^{\infty}r^{(n)}$ is convergent in \rn\ iff the diagonal series $\sum_{n=1}^{\infty}l(J^{(n)})$ is convergent in \rn\ and in this case they are equal. Thus two cases arise :

\begin{enumerate}[label=(\roman*), ref=\roman*]
\item Case $\sum_{n=1}^{\infty}r^{(n)}$ converges in \rn. Then $\sum_{n=1}^{\infty} |A_n|$ also converges in \rn, by (\ref{eq:om-count-subadd-row-n}), and we have :
\label{case:om-count-subadd-case1}
\begin{eqnarray}
|A| & \leq & \sum_{n=1}^{\infty}l(J^{(n)}), \hspace{1em} \mbox{from (\ref{eq:om-count-subadd-diag-sum})} \nonumber \\
& = & \sum_{n=1}^{\infty}r^{(n)} \leq \sum_{n=1}^{\infty} (|A_n| + \epsilon/2^n) \hspace{1em} \mbox{from (\ref{eq:om-count-subadd-row-n})} \nonumber \\
\imp |A| & \leq & \sum_{n=1}^{\infty} |A_n| + \epsilon \label{eq:om-count-subadd-rhs-epsilon}
\end{eqnarray}
\item Case $\sum_{n=1}^{\infty}r^{(n)}$ diverges in \rn. Then :
\label{case:om-count-subadd-case2}
\begin{eqnarray}
\sum_{n=1}^{\infty} (|A_n| + \epsilon/2^n) & = & \infty, \hspace{1em} \mbox{from (\ref{eq:om-count-subadd-row-n})} \nonumber \\
\imp \hspace{1em} \sum_{n=1}^{\infty} |A_n| & = & \infty \label{eq:om-count-subadd-rhs-infty}
\end{eqnarray}
\end{enumerate}

If $\exists \; \epsilon > 0$ such that case (\ref{case:om-count-subadd-case2}) arises then we obtain (\ref{eq:om-count-subadd}) immediately from (\ref{eq:om-count-subadd-rhs-infty}). Otherwise for every $\epsilon > 0$ case (\ref{case:om-count-subadd-case1}) arises, so that (\ref{eq:om-count-subadd-rhs-epsilon}) holds $\forall \; \epsilon > 0$, and hence (\ref{eq:om-count-subadd}) follows.
\end{proof}

\newpage
\subsection{Theorem (Outer Measure Unaffected By Single Point)}
\label{thm:om-add-remove-pt}
For any $A \subseteq \rn$ and $x \in \rn$ : \\

(i) $|A \cup \{x\}| = |A|$, \hspace{0.5em} (ii) $|A \sm \{x\}| = |A|$ \\

\begin{proof}
\begin{enumerate}[label=(\roman*), ref=\roman*]
\item Order Preservation (Theorem~\ref{thm:om-order-pres}) ${\imp |A| \leq |A \cup \{x\}|}$. Countable Subadditivity (Theorem~\ref{thm:om-count-subadd}) $\imp |A \cup \{x\}| \leq |A| + |\{x\}| = |A|$, since $|\{x\}| = 0$. Then by anti-symmetry of the order relation on \ern\ the result follows.
\item From (i), $|A \sm \{x\} \cup \{x\}| = |A \sm \{x\}|$, ie. $|A \cup \{x\}| = |A \sm \{x\}|$, hence $|A| = |A \sm \{x\}|.$
\end{enumerate}
\end{proof}

\subsection{Theorem (Outer Measure Equals Length for Bounded Intervals)}
\label{thm:om-len-bnd-interval}
If $I$ is a bounded interval then $|I| = l(I)$. \\

\begin{proof}
\noindent Firstly consider $I = [a, b]$ where $a < b$. The single set $(a - \epsilon, b + \epsilon)$ is a COI cover of $I$ \tf\ $|I| \leq (b - a) + 2\epsilon\ \forall \; \epsilon > 0\ \tf\ |I| \leq (b - a)$, ie $|I| \leq l(I)$. \\

Suppose now $(I_n)$ is an arbitrary COI cover for $I = [a, b]$. By Heine-Borel theorem, $I$ is compact so $\exists$ a finite subcover, hence $\exists\ m$ such that $I \subseteq I_1 \cup \cdots \cup I_m$. Then by Theorem \ref{thm:propints-finite-disj-superset} above, $l(I_1) + \cdots + l(I_m) \geq l(I)$, and so $\sum_{n=1}^{\infty} l(I_n) \geq l(I)$. Thus $l(I)$ is a lower bound for all such infinite sums, and so $|I| \geq l(I)$. Thus $|I| = l(I)$ for any $I$ of form $[a, b]$. \\

By Theorem \ref{thm:om-add-remove-pt} we can add or remove finitely many points from a set and it does not alter the outer measure. Thus for any $a, b \in \rn$ with $a < b$, we have $b - a = l([a, b]) = |[a, b]| = |(a, b)| = |[a, b)| = |(a, b]|$, thus the measure of each of the latter three interval types equals their length $b - a$. Thus $|I| = l(I)$ for every bounded interval type $I$.
\end{proof}

We can show $|I| = l(I)$ is also true for unbounded intervals, once we have shown outer measure is countably additive wrt bounded intervals (Theorem \ref{thm:om-count-add-bnd-int} below). Note that an unbounded set in \rn\ does not necessarily have outer measure of $\infty$, as outer measure is readily checked to be zero for any countable set. However due to Theorem \ref{thm:om-order-pres} and the theorem just proved, any bounded set in \rn\ must have outer measure $< \infty$.

\subsection{Theorem (Additivity for `Separated' Sets)}
\label{thm:om-separated-sets}
For any $A, B \subseteq \rn$ with $\sup A < \inf B$ (ie. a `gap' exists between $A$ and $B$), we have $|A \cup B| = |A| + |B|$. \\

\begin{proof}
\noindent By subadditivity we have $|A \cup B| \leq |A| + |B|$. For $\geq$ we require $|A| + |B|$ to be a lower bound for every $\sum_{n=1}^{\infty} l(I_n)$ where $(I_n)$ is an arbitrary COI for $A \cup B$. \\

If $\exists$ an unbounded $I_n$ then $l(I_n) = \infty$ and we are done, so we can assume every open interval $I_n$ is bounded. Furthermore if $\sum_{n=1}^{\infty} l(I_n) = \infty$ then we are again done, so assume this sum is finite. Choose any $a \in (\sup A, \inf B$).

\vspace{5ex} \input{diagram-separated-sets} \vspace{4ex}

Partition $(I_n)$ into the following subsequences : $(P_n)$ which meet $A$ only, $(R_n)$ which meet $B$ only, $(Q_n)$ which meet both $A$ and $B$ (so $a \in Q_n$), and $(T_n)$ which meet neither $A$ nor $B$. If any of these are finite subsequences extend them out to infinite sequences by adding $\es$'s (or $\{a\}$'s in the case of $(Q_n)$). \\

Then since $\sum_{n=1}^{\infty} l(I_n) < \infty$, by \cite[Theorem~1]{mse:series-rearrangements} we have :
$$
\sum_{n=1}^{\infty} l(P_n) + \sum_{n=1}^{\infty} l(Q_n) + \sum_{n=1}^{\infty} l(R_n) + \sum_{n=1}^{\infty} l(T_n) = \sum_{n=1}^{\infty} l(I_n)
$$
so
\begin{equation}
\label{eq:om-separated-sets}
\sum_{n=1}^{\infty} l(P_n) + \sum_{n=1}^{\infty} l(Q_n) + \sum_{n=1}^{\infty} l(R_n) \leq \sum_{n=1}^{\infty} l(I_n) < \infty
\end{equation}

Because every bounded open interval $Q_n$ must contain $a$ we can write $Q_n = Q_n^{(L)} \cup \{a\} \cup Q_n^{(R)}$, where $Q_n^{(L)}$ is a bounded open interval strictly to the left of $a$ and $Q_n^{(R)}$ is a bounded open interval strictly to the right of $a$. Then noting $a \notin A \cup B$ we have :
$$
A \subseteq \bigcup_{n=1}^{\infty} P_n \cup \bigcup_{n=1}^{\infty} Q_n^{(L)} \hspace{1em} \mbox{and} \hspace{1em} B \subseteq \bigcup_{n=1}^{\infty} Q_n^{(R)} \cup \bigcup_{n=1}^{\infty} R_n
$$
But then again applying \cite[Theorem~1]{mse:series-rearrangements} :
$$
|A| \leq \sum_{n=1}^{\infty} l(P_n) + \sum_{n=1}^{\infty} l(Q_n^{(L)}) \hspace{1em} \mbox{and} \hspace{1em} |B| \leq \sum_{n=1}^{\infty} l(Q_n^{(R)}) + \sum_{n=1}^{\infty} l(R_n)
$$
\tf\ since $l(Q_n) = l(Q_n^{(L)}) + l(Q_n^{(R)})\ \forall\ n$ :
\begin{eqnarray*}
|A| + |B| & \leq & \sum_{n=1}^{\infty} l(P_n) + \sum_{n=1}^{\infty} l(Q_n) + \sum_{n=1}^{\infty} l(R_n) \\
& \leq & \sum_{n=1}^{\infty} l(I_n) \hspace{2em} \mbox{by (\ref{eq:om-separated-sets})}
\end{eqnarray*}

\end{proof}

Note in the next theorem we confine ourselves to bounded intervals because it is only for these that we have presently have $|I| = l(I)$.

\subsection{Theorem (Finite Additivity wrt Bounded Intervals)}
\label{thm:om-finite-add-bnd-int}
Outer measure is finitely additive wrt bounded intervals, ie. if $I_1, \ldots, I_n$ are disjoint bounded intervals in \rn\ then :
$$
| I_1 \cup \cdots \cup I_n | = |I_1| + \cdots + |I_n|.
$$

\begin{proof}
\noindent Consider case $n = 2$. Relabel $I_1, I_2$ so they are in L to R order. If a `gap' exists between $I_1$ and $I_2$ then by Theorem \ref{thm:om-separated-sets} the result follows. If no gap exists then because $I_1$ and $I_2$ are intervals we must have $I_1 \cup I_2$ an interval with $l(I_1 \cup I_2) = l(I_1) + l(I_2)$. But since $I_1$ and $I_2$ are bounded we know $|I_1| = l(I_1)$, $|I_2| = l(I_2)$, and $|I_1 \cup I_2| = l(I_1 \cup I_2)$, from Theorem \ref{thm:om-len-bnd-interval}, and hence $|I_1 \cup I_2| = |I_1| + |I_2|$ as required. \\

Now assume the theorem is true for $n$ (where $n \geq 2$), and consider the case of $n + 1$. Let $I_1, \ldots, I_{n+1}$ be $n + 1$ disjoint bounded intervals. Relabel them so they are in L to R order. \\

If there is no `gap' between $I_n$ and $I_{n+1}$ then by adding the single boundary point $a$ between $I_n$ and $I_{n+1}$ and letting interval $J = I_n \cup \{a\} \cup I_{n+1}$ we have $l(J) = l(I_n) + l(I_{n+1})$, and :
\begin{eqnarray*}
|I_1 \cup \cdots \cup I_{n+1}| & = & |I_1 \cup \cdots \cup I_{n-1} \cup J|, \hspace{3em} \mbox{by Theorem \ref{thm:om-add-remove-pt}} \\
& = & |I_1| + \cdots + |I_{n-1}| + |J|, \hspace{1.5em} \mbox{by the inductive assumption, since } I_1, \ldots \\
& & \hspace{12em} \ldots, I_{n-1}, J  \mbox{ are $n$ disjoint bounded intervals.}
\end{eqnarray*}

But, using Theorem \ref{thm:om-len-bnd-interval} :
$$
|J| = l(J) = l(I_n) + l(I_{n+1}) = |I_n| + |I_{n+1}|
$$
hence the result for $n + 1$ follows. \\

If there is a `gap' between $I_n$ and $I_{n+1}$, ie.
$$
\sup (I_1 \cup \cdots \cup I_n) < \inf I_{n+1}
$$
then we obtain :
\begin{eqnarray*}
|I_1 \cup \cdots \cup I_{n+1}| & = & |I_1 \cup \cdots \cup I_n| + |I_{n+1}|, \hspace{2.7em} \mbox{by Theorem \ref{thm:om-separated-sets}} \\
& = & |I_1| + \cdots + |I_n| + |I_{n+1}|, \hspace{2.1em} \mbox{by the inductive assumption}
\end{eqnarray*}

hence completing the inductive step.
\end{proof}

\subsection{Theorem (Countable Additivity wrt Bounded Intervals)}
\label{thm:om-count-add-bnd-int}
Outer measure is countably additive wrt bounded intervals, ie. if $(I_n)$ is a sequence of disjoint bounded intervals in \rn\ then :
$$\leftvert{2.9ex}\bigcup_{n=1}^{\infty} I_n\rightvert{2.9ex} = \sum_{n=1}^{\infty} |I_n|$$

\begin{proof}
\noindent By countable subadditivity of outer measure (Theorem \ref{thm:om-count-subadd}) :
\begin{equation}
\label{eq:om-count-add-bnd-int}
\leftvert{2.9ex}\bigcup_{n=1}^{\infty} I_n\rightvert{2.9ex} \leq \sum_{n=1}^{\infty} |I_n|
\end{equation}

By Theorems \ref{thm:om-finite-add-bnd-int} and \ref{thm:om-order-pres} :
$$|I_1| + \cdots + |I_n| = |I_1 \cup \cdots \cup I_n| \leq \leftvert{2.9ex}\bigcup_{n=1}^{\infty} I_n\rightvert{2.9ex}, \forall\ n$$
thus, taking the limit as $n \ra\ \infty$ :
$$\sum_{n=1}^{\infty} |I_n| \leq \leftvert{2.9ex}\bigcup_{n=1}^{\infty} I_n\rightvert{2.9ex}$$
which together with (\ref{eq:om-count-add-bnd-int}) gives the required result.
\end{proof}

\subsection{Theorem (Outer Measure Equals Length for Unbounded Intervals)}
\label{thm:om-len-unbnd-interval}
If $I \subseteq \rn$ is an unbounded interval then $|I| = \infty$, ie $|I| = l(I)$. \\

\begin{proof}
\noindent First consider $I = (a, \infty)$. Then :
$$I = \bigcup_{n=1}^{\infty} (a + n - 1, a + n] = \mbox{disjoint union of bounded intervals}$$
\tf\ by Theorems \ref{thm:om-count-add-bnd-int} and \ref{thm:om-len-bnd-interval}, $|I| = \sum_{n=1}^{\infty} 1 = \infty$. \\

The case of $(-\infty, a)$ follows similarly. The cases of $[a, \infty)$ and $(-\infty, a]$ then follow by adding a single point. Finally, $I = (-\infty, \infty)$ follows from order preservation.
\end{proof}

Hence we now have $|I| = l(I),\ \forall$ intervals $I \subseteq \rn$, a property that we would expect intuitively from any measure function on \rn.

\subsection{Theorem (Countable Additivity wrt All Intervals)}
\label{thm:om-count-add-all-int}
Outer measure is countably additive wrt all intervals in \rn, ie. if $(I_n)$ is any sequence of disjoint intervals in \rn\ then :
$$\leftvert{2.9ex}\bigcup_{n=1}^{\infty} I_n\rightvert{2.9ex} = \sum_{n=1}^{\infty} |I_n|.$$

\begin{proof}
\noindent If every $I_n$ is bounded then we are done by Theorem \ref{thm:om-count-add-bnd-int}. Otherwise $\exists$ an unbounded $I_n$ and then by Theorem \ref{thm:om-len-unbnd-interval} we have $|I_n| = \infty$. Then $\sum_{n=1}^{\infty} |I_n| = \infty$. But by order preservation $\leftvert{2ex}\bigcup_{n=1}^{\infty} I_n\rightvert{2ex} = \infty$, and thus
\end{proof}

\subsection{Theorem (Countable Additivity on System of Intervals \ci)}
\label{thm:om-count-add-sys-int}
Outer measure is countably additive on the system of intervals \ci\ $= \{$ all countable unions of intervals in \rn\ $\}$. \\

\begin{proof}
\noindent Take $A = \bigcup_{n=1}^{\infty} A_n$ a disjoint countable union in \ci\ (noting \ci\ is closed under countable union). By Theorem \ref{thm:propints-disj-count-union}, each $A_n \in \ci$ is a \emph{disjoint} countable union of intervals in \rn\ : \\

\input{diagram-countable-add-1}
{ \vspace{-6ex} \begin{equation} \label{eq:om-count-add-sys-int-1} \end{equation} } 

\vspace{3ex}

We require $|A| = \sum_{n=1}^{\infty} |A_n|$, and this is clear if $\exists\ |A_n| = \infty$, so assume every $|A_n| \in \rn$. By Theorem \ref{thm:om-count-add-all-int} we have : \\

\input{diagram-countable-add-2}
{ \vspace{-6.7ex} \begin{equation} \label{eq:om-count-add-sys-int-2} \end{equation} } 

\vspace{3ex}

Label as $(J_n)$ the diagonal sequence of intervals shown in red in (\ref{eq:om-count-add-sys-int-1}), so $A = \bigcup_{n=1}^{\infty} J_n$. Then the $J_n$ are mutually disjoint intervals, and so by Theorem \ref{thm:om-count-add-all-int}, we have $|A| = \sum_{n=1}^{\infty} |J_n|$. \\

But since every $|A_n| \in \rn$, we can apply \cite[Theorem~3]{mse:series-rearrangements} to (\ref{eq:om-count-add-sys-int-2}) to conclude $\sum_{n=1}^{\infty} |J_n|$ converges in \rn\ iff $\sum_{n=1}^{\infty} |A_n|$ converges in \rn\ and in this case the two are equal. In the case they diverge in \rn\ they both equal $\infty$, thus in all cases they are equal in \ern. Thus we have the desired $|A| = \sum_{n=1}^{\infty} |A_n|$. \vspace{1ex}
\end{proof}

\subsection{The Borel Algebra}
\label{sec:borel-algebra}

In the answer \cite{mse:countable-additivity} it is discussed how the properties required of a function $\mu : \cc \ra [0, \infty]$, defined on a class \cc\ of subsets of \rn, for it to be considered a `measure' on \rn, are as follows :

\begin{enumerate}[label=(\Alph*), ref=\Alph*]
\item For the domain \cc\ of $\mu$ :
\label{list:C-props}
\begin{enumerate}[label=\arabic*., ref=\arabic*]
\item contains all intervals
\item is closed under countable union
\item is closed under set complementation
\label{list:set-compl}
\end{enumerate}
\item For the function $\mu$ :
\label{list:MU-props}
\begin{enumerate}[label=\arabic*., ref=\arabic*]
\item satisfies $\mu(I) = l(I)\ \forall$ intervals $I \subseteq \rn$
\item is countably additive, ie $\mu( \bigcup_{n=1}^{\infty} E_n ) = \sum_{n=1}^{\infty} \mu(E_n)$, whenever $(E_n)$ are disjoint sets in \cc.
\end{enumerate}
\end{enumerate}

From the theorems above, and from the Appendix~\ref{sec:sys-intervals-not-closed}, it is clear that in taking $\mu =$ outer measure and $\cc = \ci$, all of these properties are satisfied with the exception of (\ref{list:C-props})(\ref{list:set-compl}). It is shown in \cite[\S2.18, p21]{axler:measure-theory} that countable additivity of outer measure fails on $P(\rn)$ even for just two sets. Thus for outer measure to become a measure we seek a \cc\ somewhere between \ci\ and $P(\rn)$ such that all of the above conditions (\ref{list:C-props}) -- (\ref{list:MU-props}) hold. \\

Conditions (\ref{list:C-props}) give rise to the concept of a \sa\ on a set $X$, which is a non-empty collection of subsets of $X$ closed under countable union and closed under set complementation (from which closure under countable intersection follows). Then the intersection of any family of \sa s is a \sa, and if \ca\ is an arbitrary collection of subsets of $X$ we can define $\sg(\ca)$ to be the intersection of all \sa s of $X$ which contain \ca\ (noting that there is always at least one such \sa, namely $P(X)$). This is then the smallest \sa\ containing \ca. \\

If the domain \cc\ was a \sa\ containing \ci\ then (\ref{list:C-props}) would be satisfied, and $\mu : \cc \ra [0, \infty]$ would then be a measure if $\mu$ was countably additive on \cc. The Borel algebra \cb\ is defined to be the \sa\ $\sg(\ci)$. It is readily shown that $\cb\ = \sg(\ci) = \sg(\mbox{all intervals}) = \sg(\mbox{all open intervals}) = \sg(\mbox{all closed intervals}) = \sg(\mbox{all open sets}) = \sg(\mbox{all closed sets})$. In \cite[\S2D, p47]{axler:measure-theory} it is proved that the outer measure is in fact countably additive on \cb, and thus we have outer measure as a measure on \rn\ with domain \cb. The sets in \cb\ are called Borel sets. \\

Although we would naturally want the measure function to be defined on as large a class of sets as possible, taking \cb\ to be the \emph{smallest} \sa\ containing \ci\ has the advantage of providing an avenue of proving properties about the sets in \cb\ --- namely by considering the set of all subsets of \rn\ which have the desired property, and showing that it is a \sa\ and that it contains \ci\ --- then it must contain \cb\ and therefore every set in \cb\ has the desired property. An example of this technique is in \cite[\S2.65, p48]{axler:measure-theory}, where it is used to prove that the Borel sets can be approximated from below by closed sets by showing that the set \cl\ of all sets possessing that property is a \sa\ which contains all closed sets, and hence contains \cb\ $= \sg(\mbox{all closed sets})$. The sets in \cl\ are known as the Lebesgue measurable sets \cite[p52]{axler:measure-theory}.


\newpage
\input{IntervalsAndOuterMeasureOnRAppendices}


\newpage

\setcounter{secnumdepth}{0}

\section{References}

All web links retrieved on 2\textsuperscript{nd} October 2023.

\begingroup
\renewcommand{\section}[2]{}

\endgroup

\end{document}

%% file: diagram-props-intervals.tex
\begin{tikzpicture} [
	xscale=1.7, yscale=-1.7,
	every node/.style={scale=1, node font=\normalsize},
    raised term style/.style={above=-2.1ex},
	single arrow/.style=-{Stealth[length=#1 7]},
	diagonal color/.style={red},
	diagonal line/.style={line width=0.7pt, single arrow=1.0mm, diagonal color, opacity=0.35},
]


\node at (0, 0) {};
\begin{scope}[xshift=0cm]
\begin{scope}[xshift=0.5cm, yshift=0.35cm]
\node at (0.1, 0) {$A$};
\node at (0.5, 0) {$=$};
\node[raised term style] at (1, 0) {$J_1^{(1)}$};
\node at (1.5, 0) {$\cup$};
\node[raised term style] at (2, 0) {$J_2^{(1)}$};
\node at (2.5, 0) {$\cup$};
\node at (3, 0) {$\cdots$};
\node at (3.5, 0) {$\cup$};
\node[raised term style] at (4, 0) {$J_{r_1}^{(1)}$};

\node at (0.5, 0.7) {$\cup$};
\node[raised term style] at (1, 0.7) {$J_1^{(2)}$};
\node at (1.5, 0.7) {$\cup$};
\node[raised term style] at (2, 0.7) {$J_2^{(2)}$};
\node at (2.5, 0.7) {$\cup$};
\node at (3.5, 0.7) {$\cdots \cdots \cdots \cdots$};
\node at (4.5, 0.7) {$\cup$};
\node[raised term style] at (5, 0.7) {$J_{r_2}^{(2)}$};

\node at (0.5, 1.2) {$\vdots$};
\node at (2, 1.2) {$\vdots$};

\node at (0.5, 1.85) {$\cup$};
\node[raised term style] at (1, 1.85) {$J_1^{(n)}$};
\node at (1.5, 1.85) {$\cup$};
\node[raised term style] at (2, 1.85) {$J_2^{(n)}$};
\node at (2.5, 1.85) {$\cup$};
\node at (3, 1.85) {$\cdots$};
\node at (3.5, 1.85) {$\cup$};
\node[raised term style] at (4, 1.85) {$J_{r_n}^{(1)}$};

\node at (0.5, 2.4) {$\vdots$};
\node at (2, 2.4) {$\vdots$};

\end{scope}

\draw[diagonal line] ($(1, 0.7)!0.1cm!(2, 0)$) -- ($(2, 0)!0.1cm!(1, 0.7)$);
\draw[diagonal line] ($(1, 1.4)!0.1cm!(3, 0)$) -- ($(3, 0)!0.1cm!(1, 1.4)$);
\draw[diagonal line] ($(1, 2.55)!0.1cm!(3, 1.15)$) -- ($(3, 1.15)!0.8cm!(1, 2.55)$);

\end{scope}

\end{tikzpicture}

%% file: diagram-countable-subadd-1.tex
\begin{tikzpicture} [
	xscale=1.7, yscale=-1.7,
	every node/.style={scale=1, node font=\normalsize},
    raised term style/.style={above=-2.1ex},
	single arrow/.style=-{Stealth[length=#1 7]},
	diagonal color/.style={red},
	diagonal line/.style={line width=0.7pt, single arrow=1.0mm, diagonal color, opacity=0.35},
]


\node at (0, 0) {};
\begin{scope}[xshift=0.3cm]
\begin{scope}[xshift=0.5cm, yshift=0.35cm]
\node at (0, 0) {$A_1$};
\node at (0.5, 0) {$\subseteq$};
\node[raised term style] at (1, 0) {$I_1^{(1)}$};
\node at (1.5, 0) {$\cup$};
\node[raised term style] at (2, 0) {$I_2^{(1)}$};
\node at (2.5, 0) {$\cup$};
\node[raised term style] at (3, 0) {$I_3^{(1)}$};
\node at (3.5, 0) {$\cup$};
\node at (4, 0) {$\cdots$};

\node at (0, 0.7) {$A_2$};
\node at (0.5, 0.7) {$\subseteq$};
\node[raised term style] at (1, 0.7) {$I_1^{(2)}$};
\node at (1.5, 0.7) {$\cup$};
\node[raised term style] at (2, 0.7) {$I_2^{(2)}$};
\node at (2.5, 0.7) {$\cup$};
\node at (3, 0.7) {$\cdots$};

\node at (0, 1.2) {$\vdots$};
\node at (1, 1.2) {$\vdots$};

\node at (0, 1.85) {$A_n$};
\node at (0.5, 1.85) {$\subseteq$};
\node[raised term style] at (1, 1.85) {$I_1^{(n)}$};
\node at (1.5, 1.85) {$\cup$};
\node[raised term style] at (2, 1.85) {$I_2^{(n)}$};
\node at (2.5, 1.85) {$\cup$};
\node at (3, 1.85) {$\cdots$};

\node at (0, 2.4) {$\vdots$};
\node at (1, 2.4) {$\vdots$};

\end{scope}

\draw[diagonal line] ($(1, 0.7)!0.1cm!(2, 0)$) -- ($(2, 0)!0.1cm!(1, 0.7)$);
\draw[diagonal line] ($(1, 1.4)!0.1cm!(3, 0)$) -- ($(3, 0)!0.1cm!(1, 1.4)$);
\draw[diagonal line] ($(1, 2.55)!0.1cm!(3, 1.15)$) -- ($(3, 1.15)!0.8cm!(1, 2.55)$);

\end{scope}

\end{tikzpicture}

%% file: diagram-countable-subadd-2.tex
\begin{tikzpicture} [
	xscale=2, yscale=-2,
	every node/.style={scale=1, node font=\normalsize},
    raised term style/.style={above=-2.1ex},
	single arrow/.style=-{Stealth[length=#1 7]},
	diagonal color/.style={red},
	diagonal line/.style={line width=0.7pt, single arrow=1.0mm, diagonal color, opacity=0.35},
]


\node at (0, 0) {};
\begin{scope}[xshift=0.5cm]
\begin{scope}[xshift=0.5cm, yshift=0.35cm]
\node at (0, 0) {$l(I_1^{(1)})$};
\node at (0.5, 0) {$+$};
\node at (1, 0) {$l(I_2^{(1)})$};
\node at (1.5, 0) {$+$};
\node at (2.3, 0) {$\cdots\cdots\cdots$};
\node at (4.5, 0) {$= \hspace{1em} r^{(1)}$};

\node at (0, 0.7) {$l(I_1^{(2)})$};
\node at (0.5, 0.7) {$+$};
\node at (1, 0.7) {$l(I_2^{(2)})$};
\node at (1.5, 0.7) {$+$};
\node at (2.3, 0.7) {$\cdots\cdots\cdots$};
\node at (4.5, 0.7) {$= \hspace{1em} r^{(2)}$};

\node at (0, 1.2) {$\vdots$};
\node at (1, 1.2) {$\vdots$};
\node at (4.5, 1.2) {$\vdots$};

\node at (0, 1.85) {$l(I_1^{(n)})$};
\node at (0.5, 1.85) {$+$};
\node at (1, 1.85) {$l(I_2^{(n)})$};
\node at (1.5, 1.85) {$+$};
\node at (2.3, 1.85) {$\cdots\cdots\cdots$};
\node at (4.5, 1.85) {$= \hspace{1em} r^{(n)}$};

\node at (0, 2.4) {$\vdots$};
\node at (1, 2.4) {$\vdots$};
\node at (4.5, 2.4) {$\vdots$};

\end{scope}

\draw[diagonal line] ($(0, 0.7)!0.0cm!(1, 0)$) -- ($(1, 0)!0.0cm!(0, 0.7)$);
\draw[diagonal line] ($(0, 1.4)!0.0cm!(2, 0)$) -- ($(2, 0)!0.0cm!(0, 1.4)$);
\draw[diagonal line] ($(0, 2.55)!0.0cm!(1, 1.85)$) -- ($(1, 1.85)!-0.25cm!(0, 2.55)$);

\end{scope}

\end{tikzpicture}

%% file: diagram-separated-sets.tex
\begin{tikzpicture} [
	scale=1,
	every node/.style={scale=1, node font=\normalsize},
	single arrow/.style=-{Stealth[length=#1 7]},
	dashed line/.style = {dash pattern=on 2pt off 1.5pt, opacity=0.8, line width=0.4pt, red},
]

\draw (0, 0) -- (12, 0);
\filldraw (6, 0) circle[radius=0.05cm];
\filldraw (4, 0) circle[radius=0.05cm];
\filldraw (8, 0) circle[radius=0.05cm];

\node[below=8pt] at (6, 0) {$a$};
\node[below=6pt] at (4, 0) {$\sup A$};
\node[below=6pt] at (8, 0) {$\inf B$};

\draw[dashed line] (4, 0.4) -- (4, 1.6);
\draw[dashed line] (8, 0.4) -- (8, 1.6);
\draw[single arrow=1mm] (4, 1) -- (1, 1);
\draw[single arrow=1mm] (8, 1) -- (11, 1);

\node[above=8pt] at (2, 0.9) {$A$};
\node[above=8pt] at (10, 0.9) {$B$};

\node at (-0.8, 0) {};

\end{tikzpicture}

%% file: diagram-countable-add-1.tex
\begin{tikzpicture} [
	xscale=1.7, yscale=-1.7,
	every node/.style={scale=1, node font=\normalsize},
    raised term style/.style={above=-2.1ex},
	single arrow/.style=-{Stealth[length=#1 7]},
	diagonal color/.style={red},
	diagonal line/.style={line width=0.7pt, single arrow=1.0mm, diagonal color, opacity=0.35},
]


\node at (0, 0) {};
\begin{scope}[xshift=-0.3cm]
\begin{scope}[xshift=0.5cm, yshift=0.35cm]
\node at (0, 0) {$A_1$};
\node at (0.5, 0) {$=$};
\node[raised term style] at (1, 0) {$I_1^{(1)}$};
\node at (1.5, 0) {$\cup$};
\node[raised term style] at (2, 0) {$I_2^{(1)}$};
\node at (2.5, 0) {$\cup$};
\node[raised term style] at (3, 0) {$I_3^{(1)}$};
\node at (3.5, 0) {$\cup$};
\node at (4, 0) {$\cdots$};
\node at (5.5, 0) {(disjoint union)};

\node at (0, 0.7) {$A_2$};
\node at (0.5, 0.7) {$=$};
\node[raised term style] at (1, 0.7) {$I_1^{(2)}$};
\node at (1.5, 0.7) {$\cup$};
\node[raised term style] at (2, 0.7) {$I_2^{(2)}$};
\node at (2.5, 0.7) {$\cup$};
\node at (3, 0.7) {$\cdots$};
\node at (5.5, 0.7) {(disjoint union)};

\node at (0, 1.2) {$\vdots$};
\node at (1, 1.2) {$\vdots$};
\node at (5.5, 1.2) {$\vdots$};

\node at (0, 1.85) {$A_n$};
\node at (0.5, 1.85) {$=$};
\node[raised term style] at (1, 1.85) {$I_1^{(n)}$};
\node at (1.5, 1.85) {$\cup$};
\node[raised term style] at (2, 1.85) {$I_2^{(n)}$};
\node at (2.5, 1.85) {$\cup$};
\node at (3, 1.85) {$\cdots$};
\node at (5.5, 1.85) {(disjoint union)};

\node at (0, 2.4) {$\vdots$};
\node at (1, 2.4) {$\vdots$};
\node at (5.5, 2.4) {$\vdots$};

\end{scope}

\draw[diagonal line] ($(1, 0.7)!0.1cm!(2, 0)$) -- ($(2, 0)!0.1cm!(1, 0.7)$);
\draw[diagonal line] ($(1, 1.4)!0.1cm!(3, 0)$) -- ($(3, 0)!0.1cm!(1, 1.4)$);
\draw[diagonal line] ($(1, 2.55)!0.1cm!(3, 1.15)$) -- ($(3, 1.15)!0.8cm!(1, 2.55)$);

\end{scope}

\end{tikzpicture}

%% file: diagram-countable-add-2.tex
\begin{tikzpicture} [
	xscale=1.7, yscale=-1.7,
	every node/.style={scale=1, node font=\normalsize},
    raised term style/.style={above=-2.1ex},
	single arrow/.style=-{Stealth[length=#1 7]},
	diagonal color/.style={red},
	diagonal line/.style={line width=0.7pt, single arrow=1.0mm, diagonal color, opacity=0.35},
]


\node at (0, 0) {};
\begin{scope}[xshift=0.0cm]
\begin{scope}[xshift=0.5cm, yshift=0.35cm]
\node at (0, 0) {$|A_1|$};
\node at (0.5, 0) {$=$};
\node[raised term style] at (1, 0) {$|I_1^{(1)}|$};
\node at (1.5, 0) {$+$};
\node[raised term style] at (2, 0) {$|I_2^{(1)}|$};
\node at (2.5, 0) {$+$};
\node[raised term style] at (3, 0) {$|I_3^{(1)}|$};
\node at (3.5, 0) {$+$};
\node at (4, 0) {$\cdots$};

\node at (0, 0.7) {$|A_2|$};
\node at (0.5, 0.7) {$=$};
\node[raised term style] at (1, 0.7) {$|I_1^{(2)}|$};
\node at (1.5, 0.7) {$\cup$};
\node[raised term style] at (2, 0.7) {$|I_2^{(2)}|$};
\node at (2.5, 0.7) {$+$};
\node at (3, 0.7) {$\cdots$};

\node at (0, 1.2) {$\vdots$};
\node at (1, 1.2) {$\vdots$};

\node at (0, 1.85) {$|A_n|$};
\node at (0.5, 1.85) {$=$};
\node[raised term style] at (1, 1.85) {$|I_1^{(n)}|$};
\node at (1.5, 1.85) {$+$};
\node[raised term style] at (2, 1.85) {$|I_2^{(n)}|$};
\node at (2.5, 1.85) {$+$};
\node at (3, 1.85) {$\cdots$};

\node at (0, 2.4) {$\vdots$};
\node at (1, 2.4) {$\vdots$};

\end{scope}

\draw[diagonal line] ($(1, 0.7)!0.1cm!(2, 0)$) -- ($(2, 0)!0.1cm!(1, 0.7)$);
\draw[diagonal line] ($(1, 1.4)!0.1cm!(3, 0)$) -- ($(3, 0)!0.1cm!(1, 1.4)$);
\draw[diagonal line] ($(1, 2.55)!0.1cm!(3, 1.15)$) -- ($(3, 1.15)!0.8cm!(1, 2.55)$);

\end{scope}

\end{tikzpicture}

%% file: IntervalsAndOuterMeasureOnRAppendices.tex
\begin{appendices}


\section{\textcolor{white}{ - Convergence Theorems For Sequences Of Riemann Integrals}}

{\Large \bfseries Convergence Theorems For Sequences Of Riemann Integrals} \\

\label{sec:convergence-theorems}

Proofs of the Bounded Convergence, Monotone Convergence, and Dominated Convergence Theorems for Riemann integrals are given below, showing the role played by the outer measure by solving the intervals problem described in Example~3 of \S\ref{sec:need-for-outer-measure}, which allows the proof of the BCT. The MCT and DCT then follow from the BCT as corollaries. The DCT for Riemann integrals is used in the proof of Stirling's Formula in \cite{conrad:stirlings-formula}. \\

The methods of proof used are as given in \cite[\S22.13, p287]{bartle:ERA}, \cite[\S22.14, p288]{bartle:ERA}, \cite[\S22.15, p289]{bartle:ERA}, and \cite[\S25.21, p359]{bartle:ERA}, but additionally we use the equivalence of the definition of the Riemann integral in terms of partition refinement in \cite[\S22.2, p277]{bartle:ERA} (as a special case of the Riemann-Stieltjes integral) with the definition of the Riemann integral in terms of upper and lower Riemann sums in \cite[\S6.1.6, p239]{bartle:IRA} (this equivalence is established by theorems \cite[\S6.4.2, p268]{bartle:IRA} and \cite[\S6.4.3, p269]{bartle:IRA}). Note the IRA definition \cite[\S6.1.6, p239]{bartle:IRA} of Riemann integrable \emph{requires} $f$ to be bounded, whilst the ERA definition \cite[\S22.2, p277]{bartle:ERA} \emph{implies} $f$ is bounded. \\

Each of the BCT, MCT, and DCT apply to a sequence $(f_n)$ of Riemann integrable functions with pointwise limit a Riemann integrable function $f$, and have as conclusion that $\lim_{n \ra \infty} \int_a^b f_n = \int_a^b f$ (where $a = -\infty$ or $b = \infty$ or both of these, in the case of the DCT). Each theorem places a certain requirement on the $(f_n)$ :

\begin{enumerate}[label=(\roman*), ref=\roman*]
\item For BCT, $(\|f_n\|)$ has some uniform bound $B$,
\item For MCT, $(f_n)$ is a monotone sequence,
\item For DCT, $(f_n)$ is dominated by an infinitely integrable function $M$.
\end{enumerate}

\subsection{Bounded Convergence Theorem}

The first theorem below is the intervals problem of Example~3 of \S\ref{sec:need-for-outer-measure} which is proved using the outer measure.

\subsubsection[Theorem A.1.1]{Theorem (Intervals Problem)}

\label{thm:BCT-intervals-problem}

If $\delta > 0$ and if $(E_n)$ is a sequence of subsets of \rn\ with each $E_n$ containing a finite number of non-overlapping closed intervals with a total length $\geq \delta$, then there exists a point belonging to infinitely many of the $E_n$. \\

\begin{proof-sm}
\noindent We use the result in Axler [1, \S2.60, p44] on the measure of a decreasing intersection, applied to the outer measure on $\mathbb{R}$ on the Borel sets \cb. We require to prove that the set $E = \bigcap_{n=1}^{\infty} \bigcup_{m=n}^{\infty} E_m$ is non-empty, and if we denote by $F_n$ the union of the above intervals within $E_n$ (so $F_n$ is a Borel set), and let $F = \bigcap_{n=1}^{\infty} \bigcup_{m=n}^{\infty} F_m \subseteq E$, then we have $|F| = \lim_{n \ra \infty} |\bigcup_{m=n}^{\infty} F_m|$, where $\forall\ n\ |\bigcup_{m=n}^{\infty} F_m| \geq |F_n| \geq \delta$ (noting $F_n$ equals a finite union of \emph{disjoint} intervals with a total length $\geq \delta$), so that $|F| \geq \delta > 0$ and hence $F \neq \emptyset$ and thus $E \neq \emptyset$, as required. 
\end{proof-sm}

\subsubsection[Theorem A.1.2]{Theorem (\cite[Lemma 22.13, p287]{bartle:ERA})}

\label{thm:lemma-bartle}

Let $f : J \ra \rn$ be a bounded non-negative function which is Riemann integrable on $J = [0, 1]$ with $\alpha = \int_0^1 f > 0$. Then $E = f^{-1} ([\alpha / 3, \infty))$ contains a finite number of non-overlapping intervals with a total length $\geq \alpha / 3\|f\|$. \\

\begin{proof-sm}
\noindent Choose a partition $P = (a_0, \ldots, a_n)$ of $J$ such that for any intermediate tags $\frac{2}{3}\alpha < S(P, f) < \frac{4}{3}\alpha$. For each $k$, choose a tag $c_k \in [a_{k-1}, a_k]$ as follows :

\begin{enumerate}[label=(\roman*), ref=\roman*]
\item $c_k$ with $f(c_k) < \alpha / 3$, if such $c_k \in [a_{k-1}, a_k]$ exists,
\item if no such $c_k$ exists then choose any $c_k \in [a_{k-1}, a_k]$.
\end{enumerate}

In the case of (ii) we have $[a_{k-1}, a_k] \subseteq E$. If the case of (ii) never arose we would have $f(c_k) < \alpha / 3$ for $k = 1, 2, \ldots, n$ and hence $S(P, f) \leq \alpha / 3$. Hence case (ii) must arise at least once, and if $L = $ sum of lengths of the non-overlapping subintervals of $J$ in (ii) then $L > 0$. \\

With the above choice of tags $(c_k)$ we then have $S(P, f) = S_1 + S_2$, where $S_1 \leq \alpha / 3$ (from (i)), and $S_2 \leq \|f\| L$ (from (ii)). Thus :
\begin{eqnarray*}
S_1 + S_2 & > & \textstyle{\frac{2}{3}} \alpha \\
\imp \hspace{2em} \|f\| L & \geq & \textstyle{\frac{2}{3}} \alpha - S_1 \geq \textstyle{\frac{1}{3}} \alpha \\
\imp \hspace{3.5em} L & \geq & \alpha / 3 \|f\|, \hspace{5em} \mbox{as required.}
\end{eqnarray*}
\end{proof-sm}

\subsubsection[Bounded Convergence Theorem]{Theorem (Bounded Convergence Theorem, \cite[\S22.14, p288]{bartle:ERA})}

Let $(f_n)$ be a sequence of bounded Riemann integrable functions $: J \ra \rn$ with pointwise limit $f$ also a bounded Riemann integrable function on $J = [a, b]$. Then if there is a uniform bound $B$ such that $\|f_n\| \leq B\ \forall\ n$ then :
$$
\lim_{n \ra \infty} \int_a^b f_n = \int_a^b f.
$$

\begin{proof-sm}
\noindent Consider first the special case where $f = 0$, each $f_n \geq 0$, and $J = [0, 1]$. Suppose that $\int_0^1 f_n \ra 0$ does not hold. We show that in this case pointwise convergence of $(f_n)$ to $0$ fails, thus establishing the contrapositive of the result. \\

By the failure of $\int_0^1 f_n \ra 0,\; \exists\ \ep > 0$ such that for no $N$ does $n > N$ imply $|\int_0^1 f_n| < \ep$, and hence $\exists$ subsequence $(f_{n_k})$ such that $\forall\ k \in \nn, \int_0^1 f_{n_k} \geq \ep$. Define $\alpha_k = \int_0^1 f_{n_k}$ and let $E_k = f_{n_k}^{-1}([\alpha_k / 3, \infty))$. By Theorem~\ref{thm:lemma-bartle}, $E_k$ contains a finite number of non-overlapping intervals with a total length $L_k \geq \alpha_k / 3 \|f_{n_k}\|$. But $\frac{1}{\|f_{n_k}\|} \geq \frac{1}{B}$ and $\alpha_k \geq \ep$ so $L_k \geq \ep / 3 B$. Then defining $\delta = \ep/3B$ and applying Theorem \ref{thm:BCT-intervals-problem} (the intervals problem), $\exists\ x \in [0, 1]$ belonging to $E_k$ for infinitely many $k$. But then for infinitely many $k$ we have :
\begin{eqnarray*}
f_{n_k}(x) & \geq & \frac{\alpha_k}{3} \geq \frac{\ep}{3} \\
\tf \lim_{k \ra \infty} f_{n_k}(x) & \neq & 0, \hspace{2em} \mbox{as required.}
\end{eqnarray*}

Now consider the case where $f = 0$, each $f_n \geq 0$, and $J = [a, b]$. Define $g_n(x) = f_n(a + x(b - a))$ for $x \in [0, 1]$. Then $g_n \geq 0$ is Riemann integrable on $[0, 1]$ with $\|g_n\| \leq B\ \forall\ n$, and $g_n(x) \ra 0\ \forall\ x$. Then by the case just proved, we have $\int_0^1 g_n \ra 0$. But, putting $t = a + x(b - a)$, so $dt/dx = b - a$, we have 
\begin{eqnarray*}
\int_a^b f_n(t)\,dt & = & \int_0^1 f_n(a + x(b - a))(b - a)\, dx \\
& = & (b - a) \int_0^1 g_n(x) \, dx \ra 0 \mbox{ as } n \ra \infty, \hspace{2em} \mbox{as required.}
\end{eqnarray*}

Finally for the general case let $g_n(x) = |f_n(x) - f(x)|$. Then $g_n(x) \geq 0$, $g_n(x) \ra 0$, and $\|g_n\| \leq B + \|f\| = B_2\ \forall\ n$, so from the previous case $\int_a^b g_n \ra 0$. From \cite[\S6.2.13, p254]{bartle:IRA}, $f_n$ and $f$ Riemann integrable \imp\ $|f_n - f|$ Riemann integrable and :
\begin{eqnarray*}
\leftvert{2.9ex} \int_a^b f_n - \int_a^b f \rightvert{2.9ex} & \leq & \int_a^b |f_n - f| = \int_a^b g_n \ra 0 \\
\tf \int_a^b f_n & \ra & \int_a^b f, \hspace{2em} \mbox{as required.}
\end{eqnarray*}

\end{proof-sm}

\subsection[Monotone Convergence Theorem]{Monotone Convergence Theorem (\cite[\S22.15, p289]{bartle:ERA})}

\noindent Let $(f_n)$ be a sequence of bounded Riemann integrable functions $: J \ra \rn$ with pointwise limit $f$ also a bounded Riemann integrable function on $J = [a, b]$. Then if $(f_n)$ is monotone :
$$
\lim_{n \ra \infty} \int_a^b f_n = \int_a^b f.
$$

\begin{proof-sm}
\noindent Case monotone increasing. Then $f_1(x) \leq f_2(x) \leq \cdots \leq f(x) = \sup_n f_n(x)$ \tf\ $f_1(x) \leq f_n(x) \leq f(x)\ \forall\ x, n$ \tf\ $|f_n(x)| \leq \max(|f_1(x)|, |f(x)|) \leq \max(\|f_1\|, \|f\|) = B$ say, $\tf \|f_n\| \leq B$. The conclusion then follows from the BCT. \\

Case monotone decreasing. Then $f_1(x) \geq f_2(x) \geq \cdots \geq f(x) = \inf_n f_n(x)$ \tf\ $f_1(x) \geq f_n(x) \geq f(x)\ \forall\ x, n$ \tf\ $\|f_n\| \leq B$ and the conclusion follows from the BCT.
\end{proof-sm}

\subsection[Dominated Convergence Theorem]{Dominated Convergence Theorem (\cite[\S25.21, p359]{bartle:ERA})}

Suppose $(f_n)$ is a sequence of functions on $[a, \infty)$ with pointwise limit $f$, and with $f_n, f$ Riemann integrable on $[a, c]\ \forall\ c > a$. Suppose the $(f_n)$ are dominated by a function $M$ which is Riemann integrable over $[a, \infty)$, ie. $|f_n(x)| \leq M(x)\ \forall\ x \in [a, \infty)$ and $n \in \nn$. Then $f_n, f$ are Riemann integrable over $[a, \infty)$ and :
$$
\lim_{n \ra \infty} \int_a^{\infty} f_n = \int_a^{\infty} f.
$$
A similar result applies for the domains $(-\infty, a]$ and $(-\infty, \infty)$. \\

\begin{proof-sm}
(Note : we define infinite Riemann integral $\int_a^{\infty} f$ to be absolutely convergent if $\int_a^{\infty} |f|$ is convergent. In this case $\int_a^{\infty} f$ is convergent and $|\int_a^{\infty} f| \leq \int_a^{\infty} |f|$). \\

Since $|f(x)| = \lim_{n \ra \infty}|f_n(x)|$, we have $|f_n|, |f| \leq M$ on $[a, \infty)$, thus by the comparison test for infinite Riemann integrals, $\int_a^{\infty} f$ and $\int_a^{\infty} f_n$ are absolutely convergent. \\

Consider a fixed $K > a$. For each $x \in [a, K]$ and $n \in \nn, |f_n(x)| \leq M(x)$. But since $M$ is Riemann integrable on $[a, K]$ it is bounded on $[a, K]$, so $\exists\ B$ such that $|f_n(x)| \leq B\ \forall\ x \in [a, K]$ and $n \in \nn$, so that $\|f_n\|_{[a, K]} \leq B\ \forall\ n \in \nn$. Then applying the BCT we have :
\begin{equation}
\label{eq:DCT-1}
\lim_{n \ra \infty} \int_a^K f_n = \int_a^K f
\end{equation}
and this holds for any $K > a$. Now for any such $K > a$ we can write :
\begin{eqnarray}
\leftvert{2.9ex} \int_a^{\infty} f_n - \int_a^{\infty} f \rightvert{2.9ex} & = & \leftvert{2.9ex} \int_a^K f_n + \int_K^{\infty} f_n - \int_a^K f - \int_K^{\infty} f \rightvert{2.9ex} \nonumber \\
& \leq & \leftvert{2.9ex} \int_a^K f_n - \int_a^K f \rightvert{2.9ex} + \leftvert{2.9ex} \int_K^{\infty} f_n - \int_K^{\infty} f \rightvert{2.9ex} \nonumber \\
& \leq & \leftvert{2.9ex} \int_a^K f_n - \int_a^K f \rightvert{2.9ex} + \int_K^{\infty} |f_n| + \int_K^{\infty} |f|, \nonumber \\
& & \mbox{from absolute convergence of $\int_K^{\infty} f_n$ and $\int_K^{\infty} f$} \nonumber \\
& \leq & \leftvert{2.9ex} \int_a^K f_n - \int_a^K f \rightvert{2.9ex} + 2 \int_K^{\infty} M, \label{eq:DCT-2} \\
& & \mbox{by comparison with function $M$ on $[K, \infty)$}. \nonumber 
\end{eqnarray}

Now let $\ep > 0$ be arbitrary. Since $M$ is Riemann integrable on $[a, \infty)\ \exists\ K$ such that $0 \leq \int_K^{\infty} M < \ep/3$. Then using equation~(\ref{eq:DCT-1}) select $N$ such that $\forall\ n > N$, $|\int_a^K f_n - \int_a^K f| < \ep/3$. Then $\forall n > N$, from (\ref{eq:DCT-2}), we have the inequality $|\int_a^{\infty} f_n - \int_a^{\infty} f| < \ep$. But this means $\lim_{n \ra \infty} \int_a^{\infty} f_n = \int_a^{\infty} f$, as required. The case $(-\infty, a]$ is proved similarly and the case $(-\infty, \infty)$ then follows from the cases $(-\infty, 0]$ and $[0, \infty)$. \\
\end{proof-sm}


\section{\textcolor{white}{ - Proving System Of Intervals Not Closed Under Set Complementation}}

{\Large \bfseries Proving \ci\ Not Closed Under Set Complementation} \\

\label{sec:sys-intervals-not-closed}

We show how the outer measure on \rn\ can be used to prove this using the outer measure properties as developed in Axler \cite[Chap 2]{axler:measure-theory}. \\

\noindent We assume the following properties of the `Cantor Set' $C$ :
\begin{enumerate}
\item $C = [0, 1] \sm \bigcup_{n=1}^{\infty} D_n$, where $(D_n)$ is a disjoint sequence of open sets in $[0, 1]$, each one a finite disjoint union of $2^{n-1}$ `middle third' open intervals, each of length $1/3^n$.
\item $C$ is uncountable.
\end{enumerate}

\noindent From finite additivity (Theorem~\ref{thm:om-finite-add-bnd-int}) and extension of interval length (Theorem~\ref{thm:om-len-bnd-interval}) we have $|[0, 1]| = 1$ and $|D_n| = 2^{n-1}/3^n$. Then :
\begin{eqnarray}
|C| & = & |[0, 1]| - \leftvert{2.9ex}\bigcup_{n=1}^{\infty} D_n\rightvert{2.9ex}, \hspace{2em} \mbox{from Axler \cite[\S2.57(b), p42]{axler:measure-theory}} \label{eq:appendix-cantor-1} \\
& = & 1 - \sum_{n=1}^{\infty} |D_n|, \hspace{2em} \mbox{using countable additivity on \ci\ (Theorem~\ref{thm:om-count-add-sys-int})} \label{eq:appendix-cantor-2} \\
& = & 1 - \frac{1}{3} \sum_{n=1}^{\infty} {\left(\frac{2}{3}\right)}^{n-1} = 0. \nonumber 
\end{eqnarray}

Now if $C$ were a countable union of intervals, $\bigcup_{n=1}^{\infty} I_n$ say, then since $C$ is uncountable at least one of the $I_n$ must have non-zero width, $I_k$ say, for otherwise $C$ would be a countable union of singletons. But then by order preservation we would have $|C| \geq |I_k| > 0$, a contradiction. Hence $C$ can never be expressed as a countable union of intervals, ie $C \notin \ci$. We now note \ci\ is closed under `subtraction' of an interval, by Theorem~\ref{thm:propints-subtract-1-int}, since $(\bigcup I_n) \sm I = \bigcup (I_n \sm I)$. Thus if the complement ${\rn \sm \bigcup_{n=1}^{\infty} D_n}$ of the set $\bigcup_{n=1}^{\infty} D_n$ in \ci\ was also in \ci\ we would conclude by subtracting the intervals $(-\infty, 0)$ and $(1, \infty)$ from it that $C \in \ci$ --- a contradiction. Thus the set $\bigcup_{n=1}^{\infty} D_n \in \ci$ provides the required counterexample. \\

Note in equation (\ref{eq:appendix-cantor-1}) we cannot apply Theorem~\ref{thm:om-count-add-sys-int} to derive this relation since we only know $\bigcup_{n=1}^{\infty} D_n \in \ci$, but we do not have $C \in \ci$. Thus we require the fuller set of measure properties from Axler where the domain is the Borel algebra \cb\ which does contain both of these sets. However equation~(\ref{eq:appendix-cantor-2}) follows from Theorem~\ref{thm:om-count-add-sys-int}, as here we are only dealing with sets in \ci. \\

\end{appendices}